\documentclass[a4paper,12pt]{llncs}
\usepackage[mathscr]{euscript}
\usepackage{amsfonts}
\usepackage{amssymb, amsmath}
\usepackage{url}

\newtheorem{thm}{Theorem}
\newtheorem{defn}[thm]{Definition}
\newcommand{\F}{\mathbb F_q}

\begin{document}
\title{Osculating Spaces of Varieties and Linear Network Codes}
\author{Johan~P.~Hansen}
\institute{Department
of Mathematics, Aarhus University\thanks{Part of this work was done while visiting Institut de Math\'ematiques de Luminy, MARSEILLE, France.}, \email{matjph@imf.au.dk}}
\date{\today}
\maketitle
\begin{abstract}

We present a general theory to obtain good linear network codes utilizing the osculating nature of algebraic varieties. In particular, we obtain from the osculating spaces of Veronese varieties explicit families of equidimensional vector spaces, in which any pair of distinct vector spaces intersect in the same dimension.

Linear network coding transmits information in terms of a basis of a vector space and the information is received as a basis of a possible altered vector space. Ralf Koetter and Frank R. Kschischang \cite{DBLP:journals/tit/KoetterK08}  introduced a metric on the set of vector spaces and showed that a minimal distance decoder for this metric achieves correct decoding if the dimension of the intersection of  the transmitted and received vector space is sufficiently large.

The proposed osculating spaces of Veronese varieties are equidistant in the above metric. The parameters of the resulting linear network codes are determined.

\end{abstract}

\subsection*{Notation}
\begin{itemize}
\item $\F$ is the finite field with $q$ elements of characteristic $p$.
\item $\mathbb F=\overline{\mathbb F_q}$ is an algebraic closure of $\F$.
\item $R_d = \mathbb F [X_0,\dots,X_n]_d$ and $R_d(\F)= \F  [X_0,\dots,X_n]_d$ the homogenous polynomials  of degree $d$ with coefficients in $\mathbb F$ and $\F$.
\item $R=\mathbb F [X_0,\dots,X_n] = \oplus_d  R_d$ and $R(\F)=\F [X_0,\dots,X_n] = \oplus_d  R_d(\F)$
\item $\mathrm{AffCone}(Y) \subseteq \mathbb F^{M+1}$ denotes the affine cone of the subvariety
$Y \subseteq \mathbb P^{M}$ and $\mathrm{AffCone}(Y)(\F)$ its $\F$-rational points.
\item$O_{k,X,P} \subseteq \mathbb P^M$ is the embedded $k$-osculating space of a variety $X \subseteq \mathbb P^M$ at the point $P \in X$ and $O_{k,X,P}(\F)$ its $\F$-rational points, see \ref{Veronese}.
\item $\mathcal V=\sigma_d(\mathbb P^n) \subseteq \mathbb P^M$ with $M=\binom{d+n}{n}-1$ is the Veronese variety, see \ref{osculating}.

\end{itemize}
For generalities on algebraic geometry we refer to \cite{MR0463157}.
\section{Introduction}

Algebraic varieties have in general an osculating structure. By Terracini's lemma \cite{Terracini}, their embedded tangent spaces tend to be in general position. Specifically, the tangent space at a generic point $P \in \overline{Q_1Q_2}$ on the secant variety of points on some secant  is spanned by the tangent spaces at $Q_1$ and $Q_2$. In general, the secant variety of points on some secant have the expected maximal dimension and therefore the tangent spaces generically span a space of maximal dimension, see \cite{MR1234494}.

This paper suggests $k$-osculating spaces including tangent spaces of algebraic varieties  as a source for constructing  linear subspaces in general position of interest for linear network coding. The $k$-osculating spaces are presented in \ref{osculating}. 

In particular, we will present the $k$-osculating subspaces of Veronese varieties and apply them to obtain linear network codes generalizing the results in \cite{DBLP:journals/corr/abs-1207-2083}. The Veronese varieties are presented in \ref{Veronese}.

\begin{defn}\label{def} Let $X \subseteq \mathbb P^M$ be a smooth projective  variety of dimension $n$ defined over the finite field $\F$ with $q$ elements. For each positive integer $k$ we define the $k$-osculating linear network code $\mathcal C_{k,X}$. The elements of the code are the linear subspaces in $\F^{M+1}$ which are the affine cones of the $k$-osculating subspaces $O_{k,X,P}(\F)$ at $\F$-rational points $P$ on $X$, as defined in \ref{osculating}. 

Specifically

\begin{equation*}
\mathcal C_{k,X} = \{\mathrm{AffCone}(O_{k,X,P})(\F)\ \vert\ P \in X(\F)\}\ .
\end{equation*}

The number of elements in $\mathcal C_{k,X}$ is by construction $\vert X(\F)\vert$, the number of $\F$-rational points on $X$.
\end{defn}

One should remark that the elements in $\mathcal C_{k,X}$ are not necessarily
equidimensional as linear vector spaces, however, their dimension is at most $\binom{k+n}{n}$.

Applying the construction to the Veronese variety $\mathcal X_{n,d}$ presented in \ref{Veronese}, we obtain a linear network code $\mathcal C_{k,\mathcal X_{n,d}}$ and the following result, which is proved in section \ref{proof}.
\begin{thm}\label{resultat}Let $n,d$ be positive integers and consider the Veronese variety $\mathcal X_{n,d} \subseteq \mathbb P^M$, with $M=\binom{d+n}{n}-1$, defined over the finite field $\F$ as in \ref{Veronese}. 

Let $\mathcal C_{k,\mathcal X_{n,d}}$
be the associated $k$-osculating linear network code, as defined in Definition \ref{def}.

The packet length of the linear network code is $\binom{d+n}{n}$, the dimension of the ambient vector space.
The number of vector spaces in the linear network code $\mathcal C_{k,\mathcal X_{n,d}}$  is $\vert \mathbb P^n(\F)\vert= 1+q+q^2+\dots+q^n$, the number of $\F$-rational points on $\mathbb P^n$.

The vector spaces $V \in \mathcal C_{k,\mathcal X_{n,d}}$ in the linear network code are equidimensional of dimension $\binom{k+n}{n}$ as linear subspaces of the ambient  $\binom{d+n}{n}$-dimensional $\F$-vector space. 

The elements in the code are equidistant in the metric $\mathrm{dist} (V_1,V_2)$ of (\ref{dist}) of Section \ref{network}. Specifically, we have the following results.

For vector spaces $V_1, V_2 \in \mathcal C_{k,\mathcal X_{n,d}}$ with $V_1 \neq V_2$
\begin{enumerate}
\item[i)] if $2k\geq d$, then $\dim_{\F}(V_1 \cap V_2) =\binom{2k-d+n}{n}$ and 
\begin{equation*}
\mathrm{dist} (V_1,V_2)=2\  \Bigg(\binom{k+n}{n}-\binom{2k-d+n}{n}\Bigg) .
\end{equation*}
\item[ii)] if $2k\leq d$, then $\dim_{\F}(V_1 \cap V_2) =0$ and 
\begin{equation*}
\mathrm{dist} (V_1,V_2)=2 \ \binom{k+n}{n} .
\end{equation*}
\end{enumerate}
\end{thm}

\subsection{Osculating spaces}\label{osculating}
\subsubsection{Principal Parts.}
Let $X$ be a smooth variety of dimension $n$ defined over the field $K$ and let $\mathscr F$ be a locally free $\mathscr O_X$-module.
The sheaves of $k$-principal parts $\mathscr P_X^k(\mathscr F)$ are locally free and if $\mathscr L$ is of rank 1, then  $\mathscr P_X^k(\mathscr L)$ is a locally free sheaf of rank $\binom{k+n}{n}$.

There are the fundamental exact sequences
\begin{equation*}
0 \rightarrow S^k\Omega_X \otimes_{\mathscr O_X} \mathscr F \rightarrow \mathscr P_X^k(\mathscr F) \rightarrow \mathscr P_X^{k-1}(\mathscr F) \rightarrow 0\ ,
\end{equation*}
where $\Omega_X$ is the sheaf of differentials on $X$ and $S^k\Omega_X$ its $k$th symmetric power.
These sequences can be used to give a local description of the sheaf principal parts. Specifically, if $\mathscr L$ is of
 rank 1, then $\mathscr P_X^k(\mathscr L)$ is a locally free sheaf of rank $\binom{k+n}{n}$.
Assume furthermore that $X$ is affine with coordinate ring $A=K[x_1,\dots,x_n]$, then $X$ and $\mathscr L$ can be identified with $A$. Also $S^k\Omega_X$ can be identified with the forms of degree $k$ in $A[dx_1,\dots,dx_n]$ in the indeterminates $dx_1,\dots dx_n$ and $\mathscr P_X^k(\mathscr L)$ with the polynomials of total degree $\leq k$ in the indeterminates $dx_1,\dots dx_n$. For arbitrary $X$, the local picture is similar, taking local coordinates $x_1,\dots,x_n$ at the point in question replacing $A$ by the completion of the local ring at that point.

In general, for each $k$ there is a canonical morphism
\begin{equation*}
d_k: \mathscr F \rightarrow \mathscr P_X^k(\mathscr F) \ .
\end{equation*}
For $\mathscr L$ of rank 1, using local coordinates as above, $d_k$ maps an element in $A$ to its truncated Taylor series
\begin{equation*}
f=f(x_1,\dots,x_n) \mapsto \sum_{\vert \alpha \vert \leq k}\frac{1}{\vert \alpha \vert !} \frac{\partial^{\vert \alpha \vert}f}{\partial x^{\alpha}}\ ,
\end{equation*}
where $\alpha =i_1i_2\dots i_n$ and $\vert \alpha \vert =i_1+i_2+\dots+i_n$.

\subsubsection{Osculating Spaces.}
Let $X$ be a smooth variety of dimension $n$ and let  $f: X \rightarrow \mathbb  P^M$ be an immersion. For $\mathscr L=f^*\mathscr O_{\mathbb P^n}(1)$  let $\mathscr P_X^k(\mathscr L)$ denote the sheaf of principal parts of order $k$. Then $\mathscr P_X^k(\mathscr L)$ is a locally free sheaf of rank $\binom{k+n}{n}$ and there are homomorphisms
\begin{equation*}
a^k: \mathscr O_X^{M+1} \rightarrow \mathscr P_X^k(\mathscr L)\ .
\end{equation*}
For $P \in X$ the morphism $a^k(P)$ defines the $k$-osculating space $O_{k,X,P}$ to $X$ at $P$ as
\begin{equation}\label{oscdef}
O_{k,X,P}:= \mathbb P (\mathrm{Im}(a^k(P))) \subseteq \mathbb P^M\ 
\end{equation}
of projective dimension at most $\binom{k+n}{n}-1$, see \cite{MR0506323}, \cite{MR1189975} and \cite{MR1068967}. For $k=1$ the osculating space is the tangent space to $X$ at $P$.

\section{The Veronese variety}\label{Veronese}
Let $R_1=\mathbb F[X_0,\dots,X_n]_1$ be the $n+1$ dimensional vector space of linear forms in $X_0,\dots,X_n$ and let $\mathbb P^{n} =\mathbb P(R_1) $ be the associated projective $n$-space over $\mathbb F$.

For each integer $d\geq 1$, consider $R_d$ the vector space of forms of degree $d$. A basis consists of the $\binom{n+d}{d}$ monomials $X_0^{d_0} X_1^{d_1} \dots X_n^{d_n}$ with $d_0+d_1+\dots+d_n=d$.
Let $\mathbb P^M=\mathbb P(R_d)$  be the associated projective space 
of dimension $M=\binom{n+d}{d}-1$.

The $d$-uple morphism of $\mathbb P^n=\mathbb P(R_1)$  to $\mathbb P ^M = \mathbb P(R_d)$ is the morphism
\begin{eqnarray*}
\sigma_d: \mathbb P^n=\mathbb P(R_1) \rightarrow& \mathbb P ^M& = \mathbb P(R_d)\\
 L \mapsto& L^d&
\end{eqnarray*}
with image the Veronese variety \begin{equation}\label{veronese}
\mathcal X_{n,d}=\sigma_d(\mathbb P^n) = \{L^d \vert\ L \in \mathbb P(R_1)\} \subseteq \mathbb P^M.
\end{equation}

\subsection{Osculating subspaces of the Veronese variety}\label{proof}
For the Veronese variety  $\mathcal X_{n,d}$ of (\ref{veronese}), the $k$-osculating subspaces of (\ref{oscdef}) with $1\leq k<d$,  at the point $P \in \mathcal X_{n,d}$ corresponding to the 1-form $L \in R_1$,
can be described explicitly as
\begin{equation}\label{oscv}
O_{k,\mathcal X_{n,d},P} = \mathbb P(\{L^{d-k} F \vert\ F\in R_k\})= \mathbb P(R_k)\subseteq \mathbb P^M
\end{equation}
of projective dimension exactly $\binom{k+n}{n}-1$, see \cite{MR0018876}, \cite{MR1873770}, \cite{MR2319156} and \cite{MR1992894}. The osculating spaces constitute a flag of linear subspaces
\begin{equation*}
O_{1,\mathcal X_{n,d},P}\subseteq O_{2,\mathcal X_{n,d},P}\subseteq \dots \subseteq O_{d-1,\mathcal X_{n,d},P}\ .
\end{equation*}

This explicit description of the $k$-osculating spaces allows us to establish the claims in Theorem \ref{resultat}.

The associated affine cone of the $k$-osculating space in (\ref{oscv}) is 
\begin{equation}
\mathrm{AffCone}(O_{k,\mathcal X_{n,d},P})(\F) = \{L^{d-k} F \vert\ F\in R_k\}
\end{equation}
of dimension $\binom{k+n}{n}$, proving the claim on the dimension of the vector spaces in the linear network code $\mathcal C_{k,\mathcal X_{n,d}}$. 

As there is one element in $\mathcal C_{k,\mathcal X_{n,d}}$ for each $\F$-rational point on $\mathbb P^n$, it follows that the number of elements in $\mathcal C_{k,\mathcal X_{n,d}}$ is
\begin{equation*}
\vert \mathcal C_{k,\mathcal X_{n,d}}\vert = \vert \mathbb P^n(\F) \vert = 1+q+q^2+\dots+q^n\ .
\end{equation*}

Finally, let $V_1, V_2 \in \mathcal C_{k,\mathcal X_{n,d}}$ with $V_1 \neq V_2$ and 
\begin{equation*}
V_i = \{L_i^{d-k} F_i \vert\ F_i\in R_k\}
\end{equation*}

If $2k\geq d$, we have
\begin{align*}
V_1 \cap V_2=& \{L_1^{d-k} F_1 \vert\ F_1\in R_k\} \cap \{L_2^{d-k} F_2 \vert\ F_2\in R_k\}\\=& \{L_1^{d-k} L_2^{d-k} G \vert\ G\in R_{2k-d}\}\ .
\end{align*}
Otherwise the intersection is trivial, proving the claims on the dimension of the intersections and the derived distances.

\section{Linear network coding}\label{network}
In linear network, coding transmission is obtained by transmitting a number of packets into the network and each packet is regarded as a vector of length $N$ over a finite field $\F$. The packets travel the network through intermediate nodes, each forwarding $\F$-linear combinations of the packets it has available. Eventually the receiver tries to infer the originally transmitted packages from the packets that are received, see \cite{DBLP:citeseer_10.1.1.11.697} and \cite{Ho06arandom}.

All packets are vectors in $\F^N$; however, Ralf Koetter and Frank R. Kschischang \cite{DBLP:journals/tit/KoetterK08} describe a transmission model in terms of linear subspaces of $\F^N$ spanned by the packets and they define a  fixed dimension {\it code} as a nonempty subset $\mathcal C \subseteq G(n,N)(\F)$ of the Grassmannian of $n$-dimensional $\F$-linear subspaces of $\F^N$.
They endowed the Grassmannian $G(n,N)(\F)$ with the metric
\begin{equation}\label{dist}
\mathrm{dist}(V_1,V_2):=\dim_{\F}(V_1+V_2)-\dim_{\F}(V_1\cap V_2),
\end{equation}
where $V_1,V _2\in G(n,N)(\F)$.

The size of the code $\mathcal C \subseteq G(n,N)(\F)$  is denoted by $\vert \mathcal C \vert$, the minimal distance by
\begin{equation}
D(\mathcal C):= \min_{V_1,V_2 \in \mathcal C, V_1 \neq V_2} \mathrm{dist}(V_1,V_2)
\end{equation}
and $\mathcal C$ is said to be of type $[N,n,\log_q \vert \mathcal C \vert, D(\mathcal C)]$. Its normalized weight is $\lambda = \frac{n}{N}$, its rate is $R= \frac{\log_q (\vert \mathcal C \vert)}{Nn}$
and its normalized minimal distance  is $\delta = \frac{D(\mathcal C)}{2n}$.

They showed that a minimal distance decoder for this metric achieves correct decoding if the dimension of the intersection of  the transmitted and received vector-space is sufficiently large. Also they obtained  Hamming, Gilbert-Varshamov and Singleton coding bounds. 
\newcommand{\etalchar}[1]{$^{#1}$}


%
%
%
%
%
%

\end{document}